\documentclass[12pt,twoside,doublespace]{article}

\usepackage{amssymb}
\usepackage{graphicx}

\def\smskip{\par\vskip 5 pt}
\def\QED{\hfill $\Box$\smskip}
\newtheorem{theorem}{Theorem}
\newtheorem{lemma}{Lemma}
\newtheorem{proposition}{Proposition}

\begin{document}

\begin{center}

\vspace{35pt}

{\Large \bf Simplified Versions}

\vspace{5pt}

{\Large \bf  of the Conditional Gradient Method}

\vspace{5pt}

\vspace{35pt}

{\sc I.V.~Konnov\footnote{\normalsize E-mail: konn-igor@ya.ru}}

\vspace{35pt}

{\em  Department of System Analysis
and Information Technologies, \\ Kazan Federal University, ul.
Kremlevskaya, 18, Kazan 420008, Russia.}

\end{center}

\vspace{35pt}

\begin{abstract}
We suggest simple modifications of the conditional gradient
method for smooth optimization problems,
which maintain the basic convergence properties, but
reduce the implementation cost of each iteration essentially.
Namely, we propose the step-size procedure without any line-search,
and inexact solution of the direction finding subproblem. Preliminary results
of computational tests confirm efficiency of
the proposed modifications.

\end{abstract}


\newpage

%1111111111111111111111111111111111111111111111111111111111111111111111111

\section{Introduction} \label{sc:1}

Let $f : D \rightarrow \mathbb{R}$ be a function defined on some
set $D$ in the real $n$-dimensional space
$\mathbb{R}^{n}$. Then one can define the usual
optimization problem of finding the minimal value of the function
$f$ over the feasible set  $D$. For brevity, we write this problem as
\begin{equation} \label{eq:1.1}
 \min \limits _{x \in D} \to f(x),
\end{equation}
its solution set is denoted by $D^{*}$ and the optimal value of
the function by $f^{*}$, i.e.
$$
f^{*} = \inf \limits _{ x \in D} f(x).
$$
Let us first consider the well known class of smooth convex optimization problems, where
the set $D$ is supposed to be convex and closed and the
function $f$ is supposed to be convex and smooth.
This class of problems is one of the most investigated and many
iterative methods were proposed for their solution. During rather long time, the efforts
were concentrated on developing more powerful and rapidly convergent methods,
such as Newton and interior point type ones, which admit
complex transformations at each iteration, but attain high accuracy
of approximations. However, new significant areas of applications related to
data mining and processing as well as allocation decisions in information
and telecommunication networks and related systems,
where large dimensionality and inexact data together with
congestion effects and scattered necessary information force us to utilize
methods, whose iteration computation expenses and accuracy requirements are rather low,
i.e., they do not utilize matrix transformations at all.
Therefore, the well known first or even zero order methods with comparatively
slow convergence may appear very useful here.

Let us turn to the conditional gradient method (CGM for short), which is one of the oldest methods
applied to the above problem. It was first suggested in
\cite{FW56} for the case when the goal function is quadratic and
further was developed by many authors; see e.g.
\cite{LP66,DR68,PD78,Dun80}. We recall that the main idea of this
method consists in linearization of the goal function. That is,
given the current iterate $x^{k}\in D$, one finds some solution
$y^{k}$ of the problem
\begin{equation} \label{eq:1.2}
\min_{y \in D} \to \langle f'(x^{k}),y \rangle
\end{equation}
and defines $p^{k}=y^{k}-x^{k}$ as a descent direction at $x^{k}$.
Taking a suitable step-size $\lambda_{k} \in (0,1]$, one sets $x^{k+1}=x^{k}+\lambda_{k}p^{k}$
and so on.

During rather long time, this method was not considered as very efficient due to its
relatively slow convergence, but it also became very popular recently.
In fact, its auxiliary linearized problems of form (\ref{eq:1.2})
appear simpler essentially than
the quadratic ones of the most other methods. Next, it
usually yields so-called sparse approximations of a solution with few non-zero
components; see e.g. \cite{Cla10,Jag13}. These properties are very
significant for the new applications indicated above. We observe that
many efforts were directed to enhance the usual
(CGM); see e.g. \cite{GM86,BT04,Jag13,FG16,Kon17} and the references therein.
The most popular way consists in developing versions that attain more rapid convergence.
At the same time, we can create more efficient methods via reduction of the
implementation costs at each iteration with preserving all the useful properties
of the initial method.

In this paper, we will follow the second way. First of all,
being based on the approach in \cite{Kon18},
we suggest a new step-size procedure in (CGM)
without any line-search. Our new step-size procedure admits
different changes of the step-size and wide variety of implementation rules,
not only decrease as in \cite{Kon18}. It does not utilize
a priori information such as Lipschitz constants, but takes into account
 behavior of the iteration sequence, unlike the well known
divergent series rule. Moreover, the Lipschitz continuity of the
gradient of the goal function is not necessary for its convergence.
Afterwards, we introduce special threshold control and tolerances in order to avoid
exact solution of the direction finding subproblem (\ref{eq:1.2}).
We establish a complexity estimate for this method, which appears equivalent to
the convergence rate of the custom (CGM).
Furthermore, we propose a version that combines both the modifications
and show that it possesses strengthened convergence properties with respect
to the first modification of (CGM) without line-search.
Preliminary results of computational experiments confirm efficiency of all
the proposed modifications.

%22222222222222222222222222222222222222222222222222222222222222222222222222222

\section{Properties of the usual conditional gradient method}\label{sc:2}

We will take the following set of basic assumptions for
 problem (\ref{eq:1.1}).

{\bf (H)} {\em $D$ is a nonempty, convex, closed, and bounded subset
of $\mathbb{R}^{n}$, $f : \mathbb{R}^{n} \rightarrow \mathbb{R}$ is a smooth function
on the set $D$.}

Together with problem (\ref{eq:1.1}), we will consider
the following variational inequality (VI for short): Find a point
$x^{*} \in D$ such that
\begin{equation} \label{eq:2.1}
\langle f'(x^{*}),x-x^{*} \rangle \geq 0 \quad \forall x \in D.
\end{equation}
We denote by $D^{0}$ the solution set of VI (\ref{eq:2.1}).

We recall that a differentiable function $\varphi : \mathbb{R}^{n}
\to \mathbb{R}$ is called {\em pseudo-convex} on a set $D \subseteq \mathbb{R}^{n}$ if for each pair of points $
x, y \in D$  we have
$$
    \langle  \varphi' (x), y - x \rangle \geq 0 \
      \Longrightarrow \
    \varphi (y) \geq \varphi (x).
$$
It is well known that that the class of convex functions is strictly contained
in the class of pseudo-convex functions. For instance, the function $\ln t$
is concave and pseudo-convex on $\mathbb{R}_{>}=\{ t \ | \ t>0 \}$,
 but it is clearly non-convex. VI (\ref{eq:2.1}) can be used as an
optimality condition for problem (\ref{eq:1.1}) so that
solutions of VI (\ref{eq:2.1}) are called {\em stationary points} of (\ref{eq:1.1}).

%----------------------lm 2.1---------------------------------

\begin{lemma}\label{lm:2.1}\cite[Theorems 5.5 and 9.12]{Kon13} Let {\bf (H)} hold.

(i)  Each solution of problem (\ref{eq:1.1}) solves VI (\ref{eq:2.1}).

(ii)  If $f$ is pseudo-convex, then  each solution of VI (\ref{eq:2.1})
solves  problem (\ref{eq:1.1}).
\end{lemma}

The boundedness of $D$ guarantees that problem (\ref{eq:1.1}) has a solution,
moreover, $D^{*}$ is a compact set, which will be also convex if $f$ is pseudo-convex.
Given a point $x \in D$ we define the auxiliary problem
\begin{equation} \label{eq:2.2}
 \min \limits _{y \in D} \to \langle  f'(x), y\rangle .
\end{equation}
We denote by $Z(x)$ the solution set of problem  (\ref{eq:2.2}),
thus defining the set-valued  mapping $x\mapsto Z(x)$. Observe that
the set $Z(x)$ is always non-empty, convex, and compact. Also, let
$$
\mu (x)=\max \limits _{y \in D}  \langle  f'(x), x-y\rangle .
$$

Given a set $V \subseteq \mathbb{R}^{n}$,
a set-valued mapping $u\mapsto Q(u)$ is said to be {\em closed} on a set
$V$, if for each pair of sequences $\{u^{k} \} \to u$,   $\{ q^{k} \} \to q$
such that $u^{k} \in V$ and $q^{k} \in Q(u^{k})$,
we have $q \in Q(u)$.

%----------------------lm 2.2---------------------------------

\begin{lemma}\label{lm:2.2} \cite[Lemma 6.3]{Kon13}
Let {\bf (H)} hold. Then the mapping  $x\mapsto Z(x)$ is closed on $D$.
\end{lemma}

%----------------------lm 2.3---------------------------------

\begin{lemma}\label{lm:2.3}\cite[Lemma 6.4]{Kon13}
Let {\bf (H)} hold. Then the following assertions are equivalent:
\begin{description}
  \item[\rm (i)] $x^* \in D^{0}$;
  \item[\rm (ii)] $x^* \in Z(x^*)$;
  \item[\rm (iii)] $\langle f'(x^*), x^*-z^* \rangle =0$ for $z^* \in
Z(x^*)$ and $x^* \in D$.
\end{description}
\end{lemma}

The above properties are very useful for substantiation of various
(CGM) type methods. Following \cite{PD78}, we now describe
the usual (CGM) with the Armijo step-size rule for a more clear comparison with the
new methods. Here and below, $\mathbb{Z}_{+}$ denotes the set of non-negative integers.

%===========================Method (CGM)==================================
\medskip
\noindent {\bf Method (CGM).}

{\em Step 0:} Choose a point
$x^{0}\in D$, numbers $\beta \in (0,1)$ and $\theta \in (0,1)$. Set $k=0$.

{\em Step 1:} Find a point $y^{k}\in Z(x^k)$, set $d^{k}=y^{k}-x^{k}$.
If $\langle f'(x^{k}), d^{k} \rangle =0$, stop.

{\em Step 2:}  Determine $m$ as the
smallest number in $\mathbb{Z}_{+}$ such that
\begin{equation} \label{eq:2.3}
  f(x^{k}+\theta ^{m} d^{k}) \leq f(x^{k})+\beta \theta ^{m}\langle f'(x^{k}), d^{k}\rangle,
\end{equation}
set $\lambda_{k}=\theta ^{m}$, $x^{k+1}=x^{k}+\lambda_{k}d^{k}$,
$k=k+1$, and go to Step 1.
\medskip

Clearly, termination of the method yields a point of $D^{0}$. For this reason,
we will consider only the non-trivial case where the sequence $\{x^{k}\}$ is infinite.
We give the basic convergence result of the above
method.

%2=========================pro 2.1============================================

\begin{proposition} \label{pro:2.1} (e.g. \cite[Theorems 6.12 and 9.12]{Kon13})
Let {\bf (H)} hold, the sequence $\{x^k \}$ be generated by (CGM).
Then:

(i) The linesearch procedure in Step 2 is always finite.

(ii) The sequence $\{x^k \}$ has limit points, all these limit points belong to the set
$D^{0}$.

(iii) If $f$ is pseudo-convex, then all the limit points of the sequence $\{x^k \}$
belong to the set $D^{*}$, besides,
\begin{equation}\label{eq:2.4}
\lim \limits_{k\rightarrow \infty} f(x^{k})=f^{*}.
\end{equation}
\end{proposition}

We can in principle take the exact one-dimensional minimization rule instead of the current
Armijo rule in (\ref{eq:2.3}), but it is not so suitable for implementation.
Next, if the gradient of the function $f$ is Lipschitz continuous on $D$ with
some constant $L >0$, i.e., $\| f'(y)-f'(x)\| \leq L \|y-x \|$
for any vectors $x$ and $y$,  one can give bounds for the step-size and obtain the
convergence rate.

%2=========================pro 2.2============================================

\begin{proposition} \label{pro:2.2} (\cite[Theorem 6.1]{LP66}
and \cite[Chapter III, Theorem 1.7]{DR68}) Suppose that
the assumptions in {\bf (H)} are fulfilled, the function $f$ is convex,
the gradient of the function $f$ is Lipschitz continuous on $D$ with
some constant $L >0$, a sequence $\{x^{k}\}$ is generated by (CGM) where
the step-size $\lambda_{k}$ is chosen by the formula
$$
\lambda_{k}=\min\{1,\theta_{k}\sigma_{k}\},
\sigma_{k}=-\langle f'(x^{k}), d^{k}\rangle/\|d^{k}\|^{2},  \ \theta_{k}\in [\theta',\theta''], \ \theta'
>0, \theta''<2/L.
$$
Then these exists some constant $C<+\infty$ such that
\begin{equation} \label{eq:2.5}
f(x^k)-f^{*}\leq C/k \quad \mbox{for} \ k=0,1,\ldots
\end{equation}
\end{proposition}

This version reduces the computational expenses essentially due to the absence of the line-search,
but requires the evaluation of the Lipschitz constant. However, utilization of
its inexact estimates usually leads to slow convergence. This is also the case for the
known divergent series rule (see e.g. \cite{DH78})
$$
 \sum \limits_{k=0}^{\infty }\lambda_{k}=\infty, \
 \sum \limits_{k=0}^{\infty }\lambda^{2}_{k}<\infty, \ \lambda_{k} \in (0, 1), \ k=0,1,2, \ldots,
$$
and for similar rules, which do not evaluate the information about the problem
along the current iterates.

%33333333333333333333333333333333333333333333333333333333333333333333333333

\section{A simple adaptive step-size without line-search} \label{sc:3}

We now describe a modification of the (CGM), which involves a
simple adaptive step-size procedure without line-search.
Moreover, it does not require any a priori information about the problem.

%===========================Method (CGMS)==================================
\medskip
\noindent {\bf Method (CGMS).}

{\em Step 0:} Choose a point
$x^{0}\in D$, numbers $\beta \in (0,1)$ and a sequence
$\{\tau _{l}\} \to 0$, $\tau _{0} \in (0, 1)$. Set $k=0$, $l=0$,
choose a number $\lambda_{0} \in (0, \tau _{0}]$.

{\em Step 1:} Find a point $y^{k}\in Z(x^k)$, set $d^{k}=y^{k}-x^{k}$.
If $\langle f'(x^{k}), d^{k} \rangle =0$, stop.

{\em Step 2:} Set $x^{k+1}=x^{k}+\lambda_{k}d^{k}$. If
\begin{equation} \label{eq:3.1}
  f(x^{k+1}) \leq f(x^{k})+\beta \lambda_{k}\langle f'(x^{k}), d^{k}\rangle,
\end{equation}
take $\lambda_{k+1} \in [\lambda_{k}, \tau _{l}]$. Otherwise
set $\lambda'_{k+1} = \min \{\lambda_{k}, \tau _{l+1}\}$,
$l=l+1$ and take $\lambda_{k+1} \in (0,\lambda'_{k+1}]$.
Set $k=k+1$ and go to Step 1.
\medskip

Again, termination of the method yields a point of $D^{0}$ due to Lemma \ref{lm:2.3}. Hence,
we will consider only the case where the sequence $\{x^{k}\}$ is infinite.

%3=========================thm 3.1============================================

\begin{theorem} \label{thm:3.1}
Let the assumptions in {\bf (H)} be fulfilled. Then:

(i) The sequence $\{x^k \}$ has a limit point, which belongs to the set
$D^{0}$.

(ii)  If $f$ is pseudo-convex, then all the limit points of the sequence $\{x^k \}$
belong to the set $D^{*}$, besides, (\ref{eq:2.4}) holds.
\end{theorem}
{\bf Proof.} First we note that both the sequences
$\{x^{k}\}$ and $\{y^{k}\}$ belong to the bounded
set $D$ and hence have limit points. Let us consider two possible
cases.

\textit{Case 1: The number of changes of the index $l$ is finite.} \\
Then we have $\lambda _{k} \geq \bar \lambda>0$ for $k$ large enough, hence
(\ref{eq:3.1}) gives
$$
 f(x^{k+1}) \leq f(x^{k})+\beta \bar \lambda \langle f'(x^{k}), d^{k}\rangle
$$
for $k$ large enough. Since $f(x^{k}) \geq f^{*}> -\infty$, we must have
\begin{equation} \label{eq:3.2}
\lim \limits_{k\rightarrow \infty }f(x^{k})=\mu
\end{equation}
and
\begin{equation} \label{eq:3.3}
\lim \limits_{k\rightarrow \infty}\langle f '(x^{k}), d^{k}\rangle=0.
\end{equation}

Let $x'$ be an arbitrary limit point of the sequence  $\{x^{k}\}$. Taking a subsequence if
necessary we have the corresponding limit point $y'$ of the sequence
 $\{y^{k}\}$, i.\,e.
$$
\lim \limits_{s\rightarrow \infty }x^{k_{s}}=x' \ \mbox{and} \ \lim
\limits_{s\rightarrow \infty }y^{k_{s}}=y'.
$$
From (\ref{eq:3.3}) we now have
$$
\langle f '(x'), y'-x'\rangle = 0,
$$
but  the mapping  $x\mapsto Z(x)$ is closed due to Lemma
\ref{lm:2.2}, hence $y' \in Z(x')$. From Lemma \ref{lm:2.3} it
follows that $x' \in D^{0}$.  Hence,
in this case all the limit points of the sequence $\{x^k \}$
belong to the set $D^{0}$. Therefore,
assertion (i) is true. If $f$ is pseudo-convex, then $D^{0}=D^{*}$
due to Lemma \ref{lm:2.1}, which gives $\mu=f^{*}$
in (\ref{eq:3.2}) and (\ref{eq:2.4}). We conclude that
assertion (ii) is  also true.

\textit{Case 2: The number of changes of the index $l$ is infinite.} \\
Then there exists an infinite subsequence of indices $\{k_{l}\}$
 such that
$$
f (x^{k_{l}}+\lambda _{k_{l}}d^{k_{l}})-f (x^{k_{l}})=f
(x^{k_{l}+1})-f (x^{k_{l}}) > \beta \lambda _{k_{l}}\langle f
'(x^{k_{l}}), d^{k_{l}}\rangle,
$$
or equivalently,
\begin{equation} \label{eq:3.4}
\frac{f (x^{k_{l}}+\lambda _{k_{l}} d^{k_{l}})-f
(x^{k_{l}})}{\lambda _{k_{l}}}
> \beta \langle f '(x^{k_{l}}), d^{k_{l}}\rangle;
\end{equation}
besides,
$$
 \lambda_{k_{l}} \in (0,\tau_{l}], \ \lambda_{k_{l}+1} \in (0,\tau_{l+1}],
$$
and
$$
\lim \limits_{s\rightarrow \infty}\tau_{l}=0.
$$

Let $\bar x$ be an arbitrary limit point  of this subsequence
$\{x^{k_{l}}\}$. Taking a subsequence if necessary we can choose the
corresponding limit point $\bar y$ of the subsequence
 $\{y^{k_{l}}\}$. Without loss of generality we can suppose that
$$
\lim \limits_{l\rightarrow \infty }x^{k_{l}}=\bar x \ \mbox{and} \ \lim
\limits_{l\rightarrow \infty }y^{k_{l}}=\bar y.
$$
Since $ \lambda _{k_{l}} \to 0$ as $l\rightarrow \infty$,
taking the limit $l\rightarrow +\infty $ in relation (\ref{eq:3.4})
we obtain
$$
\langle f '(\bar x), \bar y-\bar x\rangle \geq \beta \langle f '(\bar x),
\bar y-\bar x\rangle ,
$$
i.e.
$$
\langle f '(\bar x), \bar y-\bar x\rangle \geq 0.
$$
By Lemma \ref{lm:2.2} we have $\bar y \in Z(\bar x)$, but from Lemma
\ref{lm:2.3} it now follows that $\bar x \in D^{0}$. Therefore,
assertion (i) is true. If $f$ is pseudo-convex, then $D^{0}=D^{*}$
due to Lemma \ref{lm:2.1}. It follows that all these limit
points of the subsequence $\{x^{k_{l}}\}$ belong to the set $D^{*}$.
Since $x^{k_{l}+1}=x^{k_{l}}+\lambda _{k_{l}}d^{k_{l}}$, $ \lambda
_{k_{l}} \to 0$, and the sequence  $\{d^{k_{l}}\}$ is bounded,  the
limit points of the subsequences $\{x^{k_{l}}\}$ and
$\{x^{k_{l}+1}\}$ coincide and  all they belong to the set $D^{*}$.

For any index $k$ we define the index $m(k)$ as follows:
$$
 \ m(k)= \max\{j \ | \ j \leq k, \ f (x^{j})-f (x^{j-1}) > \beta
\lambda _{j-1}\langle f'(x^{j-1}), d^{j-1}\rangle \},
$$
i.\,e. $j$ is the closest to $k$ but not greater index from the
subsequence $\{x^{k_{l}+1}\}$. This means that $j=k$ if $f (x^{k})-f
(x^{k-1}) > \beta \lambda _{k-1}\langle f'(x^{k-1}), d^{k-1}\rangle$.
By definition, we have
\begin{equation} \label{eq:3.5}
f(x^{k})\leq f (x^{m(k)}).
\end{equation}
Let now $x'$ be an arbitrary
limit point  of the sequence $\{x^{k}\}$, i.e. $\lim
\limits_{s\rightarrow \infty }x^{t_{s}}=x'$. Create the corresponding infinite
subsequence $\{x^{m(t_{s})}\}$. From (\ref{eq:3.5}) we have $f^{*} \leq f
(x^{t_{s}})\leq f (x^{m(t_{s})})$, but all the limit points of the
sequence $\{x^{m(t_{s})}\}$ belong to the set $D^{*}$ since it is
contained in the sequence $\{x^{k_{l}+1}\}$. Choose any limit point
$x''$ of $\{x^{m(t_{s})}\}$.  Then, taking a subsequence if necessary we
obtain
$$
f^{*} \leq f (x')\leq f (x'') = f^{*}.
$$
therefore $x' \in D^{*}$. This means that all the limit points of the
sequence $\{x^{k}\}$ belong to the set  $D^{*}$ and that
(\ref{eq:2.4}) holds true. We conclude that
assertion (ii) is  also true.
\QED

It should be observed that (CGMS) follows the approach in \cite{Kon18},
but the step-size procedure in (CGMS) admits wide variety of implementation rules
in comparison with those in \cite{Kon18}, where only the strict decrease
is indicated for possible changes of the step-size.
Even the simplest implementation rule of (CGMS), where
$\lambda_{k+1} =\max \{\lambda_{k}, \tau _{l}\}$ if (\ref{eq:3.1}) holds
and $\lambda_{k+1} =\min \{\lambda_{k}, \tau _{l+1}\}$
otherwise, admits the increase of $\lambda _{k+1} $, which prevents from the
too small step-size. Such a modification seems especially significant
for the case where the computation of the goal function value is rather expensive.

%444444444444444444444444444444444444444444444444444444444444444444444444

\section{Inexact solution of the direction finding subproblem} \label{sc:4}

It was noticed in Section \ref{sc:1} that the auxiliary direction finding subproblem (\ref{eq:1.2})
in (CGM) is simpler essentially than the quadratic ones in the projection based methods.
Nevertheless, its exact solution may also be expensive. 
If the feasible set $D$ is a general polyhedron with many vertices, one has to apply
a special algorithm at each iteration. Then, the method with approximate solution of
subproblem (\ref{eq:1.2}) may appear more efficient.
There exist several versions of such methods; see e.g. \cite{DR68,DH78}.
We observe that all these versions involve evaluation of the accuracy of
a solution of subproblem (\ref{eq:1.2}), which must tend to zero.
In this section, we intend to present some other version of this modification of (CGM),
which is based on inserting tolerances and some threshold control
of the descent property. We observe that this approach was first suggested
for the bi-coordinate descent method in \cite{Kon16b}. In \cite{Kon17a}, it was applied in a
generalized conditional gradient method for optimization problems on
Cartesian product sets, where the corresponding partial
auxiliary problems in subspaces are still to be solved exactly.
We now describe the general inexact (CGM) with the same Armijo step-size rule.

%===========================Method (CGMI)==================================
\medskip
\noindent {\bf Method (CGMI).} \\
 {\em Initialization:} Choose a point $w^{0} \in D$, numbers $\beta \in (0,1)$,
$\theta  \in (0,1)$, and a positive sequence $\{\delta _{p}\} \to 0$. Set $p=1$.

{\em Step 0:} Set  $k=0$, $x^{0}=w^{p-1}$.

{\em Step 1:}  Find a point $z^{k}\in D$ such that
\begin{equation} \label{eq:4.1}
\langle f'(x^{k}), x^{k}-z^{k} \rangle  \geq \delta _{p}.
\end{equation}
If  $\mu(x^{k})<\delta _{p} $, set $w^{p}=x^{k}$, $p=p+1$ and go to Step 0. {\em (Restart)}

{\em Step 2:}  Set $d^{k}= z^{k}-x^{k}$,
determine $m$ as the smallest number in $\mathbb{Z}_{+}$ such that
\begin{equation} \label{eq:4.2}
 f (x^{k}+\theta ^{m} d^{k}) \leq f (x^{k})+\beta \theta ^{m}
  \langle f'(x^{k}),d^{k} \rangle,
\end{equation}
set $\lambda_{k}=\theta ^{m}$, $x^{k+1}=x^{k}+\lambda_{k}d^{k}$, $k=k+1$ and go to Step 1.
\medskip

Thus, the method has a two-level structure where each outer iteration (stage) $p$
contains some number of inner iterations in $k$
with the fixed tolerance $\delta _{p}$. Completing each stage,
that is marked as restart, leads to decrease of its value.
Observe that only the restart situation requires the exact solution
of the auxiliary subproblem (\ref{eq:1.2}). In all the other cases,
we can take $z^{k}\in D$ as an arbitrary suitable point
(say a vertex of $D$) within condition (\ref{eq:4.1}).

By (\ref{eq:4.1}), we have
$$
\langle f'(x^{k}),d^{k}\rangle \leq -\delta_{p}<0
$$
in (\ref{eq:4.2}). It follows that
 \begin{equation} \label{eq:4.3}
 f (x^{k+1})-f (x^{k}) \leq \beta \lambda_{k} \langle f'(x^{k}),d^{k}\rangle \leq -\beta
 \lambda_{k}\delta_{p}.
\end{equation}

We first justify the linesearch.

%============================lm 4.1=====================================================

\begin{lemma} \label{lm:4.1} Let the assumptions in {\bf (H)} be fulfilled.
Then the linesearch procedure in Step 2 of (CGMI) is always finite.
\end{lemma}
{\bf Proof.}
If we suppose that the linesearch procedure is infinite, then (\ref{eq:4.2}) does not hold and
$$
\theta ^{-m}(f (x^{k}+\theta ^{m}d^{k}) - f (x^{k}))>\beta \langle f'(x^{k}),d^{k}\rangle,
$$
for $m \to \infty$. Hence, by taking the limit we have
$    \langle f'(x^{k}),d^{k}\rangle \geq \beta \langle f'(x^{k}),d^{k}\rangle$, hence
$\langle f'(x^{k}),d^{k}\rangle \geq 0$, a contradiction
with (\ref{eq:4.1}). \QED

We show that each stage is well defined.

%=========================pro:4.1==================================================

\begin{proposition} \label{pro:4.1} Let the assumptions in {\bf (H)} be fulfilled.
Then the number of iterations at each stage $p$ is finite.
\end{proposition}
{\bf Proof.}
Fix any $p$ and suppose that the sequence $\{x^{k}\}$ is infinite.
By (\ref{eq:4.3}), we have $f^{*}\leq f(x^{k})$ and $f(x^{k+1})\leq f(x^{k})-\beta \delta_{p}
\lambda_{k}$, hence
$$
\lim \limits_{k\rightarrow \infty }\lambda_{k}=0.
$$
Both the sequences $\{x^{k}\}$ and $\{z^{k}\}$ belong to the bounded
set $D$ and hence have limit points. Without loss of generality, we can
suppose that the subsequence $\{x^{k_{s}}\}$ converges to a point
$\bar x$ and the corresponding subsequence $\{z^{k_{s}}\}$ converges to a point
$\bar z$. Due to (\ref{eq:4.1}) we have
\begin{equation} \label{eq:4.4}
\langle f'(\bar x),\bar y -\bar x \rangle= \lim \limits_{s\rightarrow \infty}
\langle f'(x^{k_{s}}), y^{k_{s}}-x^{k_{s}} \rangle \leq -\delta_{p}.
\end{equation}
However,  (\ref{eq:4.2}) does not hold for the step-size $\lambda_{k}/\theta$.
Setting $k=k_{s}$ gives
$$
(\lambda_{k_{s}}/\theta)^{-1 }(f
(x^{k_{s}}+(\lambda_{k_{s}}/\theta)  d^{k_{s}}) - f
(x^{k_{s}}))>\beta \langle f'(x^{k_{s}}), d^{k_{s}} \rangle,
$$
hence, by taking the limit $s\rightarrow \infty$ we obtain
\begin{eqnarray*}
\langle f'(\bar x),\bar y -\bar x \rangle &=& \lim \limits_{s\rightarrow \infty
} \left\{(\lambda_{k_{s}}/\theta)^{-1 }(f
(x^{k_{s}}+(\lambda_{k_{s}}/\theta) d^{k_{s}}) - f
(x^{k_{s}}))\right\} \\
&\geq&  \beta \langle f'(\bar x),\bar y -\bar x
\rangle,
\end{eqnarray*}
i.e.,  $   (1-\beta ) \langle f'(\bar x),\bar y -\bar x \rangle \geq 0$,
which is a contradiction with (\ref{eq:4.4}). \QED

We are ready to prove convergence of the whole method.

%4===================thm 4.1=========================================

\begin{theorem} \label{thm:4.1}
Let the assumptions in {\bf (H)} be fulfilled. Then:

(i)  The number of changes of index $k$ at each stage $p$  is finite.

(ii) The sequence $\{w^{p}\}$ generated by method (CGMI) has limit points, all
these limit points belong to $D^{0}$.

(iii)  If $f$ is pseudo-convex, then all the limit points of the sequence $\{w^{p}\}$
belong to the set $D^{*}$, besides,
\begin{equation} \label{eq:4.4a}
 \lim \limits_{p\rightarrow \infty} f (w^{p})=f^{*};
\end{equation}
\end{theorem}
{\bf Proof.} Assertion (i) has been obtained in Proposition \ref{pro:4.1}.
By construction, the sequence $\{w^{p}\}$ is bounded, hence it has limit points.
Moreover, $f^{*}\leq f(w^{p+1})\leq f(w^{p})$, hence
\begin{equation} \label{eq:4.5}
\lim \limits_{p\rightarrow \infty }f(w^{p})=\mu.
\end{equation}
For each $p$ and any point $u^{p}\in Z(w^{p})$ it holds that
\begin{equation} \label{eq:4.6}
\langle f'(w^{p}),w^{p}-u^{p}\rangle \leq \delta_{p}.
\end{equation}
Fix this sequence $\{u^{p}\}$. It is also bounded and must have limit points.
Take an arbitrary limit point $\bar w$ of $\{w^{p}\}$. Then,
without loss of generality we can suppose that
$$
\bar u =\lim \limits_{t\rightarrow \infty }u^{p_{t}} \ \mbox{and} \ \bar w =\lim \limits_{t\rightarrow \infty }w^{p_{t}},
$$
for some subsequences $\{u^{p_{t}}\}$ and $\{w^{p_{t}}\}$.
Taking the limit $t\rightarrow \infty$ in (\ref{eq:4.6}) with $p=p_{t}$, we obtain
$$
\langle f'(\bar w),  \bar w-\bar u \rangle =\lim \limits_{t\rightarrow \infty }
\langle f'(w^{p_{t}}),w^{p_{t}}-u^{p_{t}}\rangle \leq 0.
$$
By Lemma \ref{lm:2.2} we have $\bar u \in Z(\bar w)$, hence
$\langle f'(\bar w),  \bar w-\bar u \rangle =0$ and $\bar w \in D^{0}$
due to Lemma \ref{lm:2.3}. This means that
all the limit points of $\{u^{p}\}$ belong to $D^{0}$.
This gives assertion (ii).  If $f$ is pseudo-convex, then $D^{0}=D^{*}$
due to Lemma \ref{lm:2.1}, which gives $\mu=f^{*}$
in (\ref{eq:4.5}) and (\ref{eq:4.4a}). We conclude that
assertion (iii) is true.
\QED

It was observed in Section \ref{sc:2} that the usual (CGM) attains
the convergence rate $O(1/k)$ under the additional
assumptions that the function $f$ is convex and its
gradient is Lipschitz continuous; see
Proposition \ref{pro:2.2} and formula (\ref{eq:2.5}). This means that
the total number of iterations $N(\varepsilon)$ that is necessary for
attaining some prescribed
accuracy  $\varepsilon>0$  for the gap value
$\Delta(x)=f(x)-f^{*}$ is estimated as follows:
\begin{equation} \label{eq:4.7}
N(\varepsilon)\leq C/\varepsilon, \ \mbox{where} \ 0<C<\infty.
\end{equation}
We can try to obtain a similar estimate for (CGMI) with the proper specialization.
In fact, if the gradient of the function $f$ is Lipschitz continuous on $D$ with
some constant $L >0$,  we can take the well known property of
such functions
$$
f (y) \leq f(x)+\langle
f'(x),y-x \rangle +0.5L\|y-x \|^{2};
$$
see \cite[Chapter III, Lemma 1.2]{DR68}. Then,
at Step 2 we have
$$
 f (x^{k}+\lambda d^{k}) - f (x^{k}) \leq \lambda [\langle f'(x^{k}),d^{k} \rangle
    +0.5L \lambda\|d^{k}\|^{2} ]  \leq \beta \lambda \langle f'(x^{k}),d^{k} \rangle,
$$
if $ \lambda \leq -2(1-\beta)\langle f'(x^{k}),d^{k} \rangle /(L\|d^{k}\|^{2})$.
Next, $\langle f'(x^{k}),d^{k} \rangle \leq -\delta_{p}$ at stage $p$,
besides, $ \|d^{k}\| \leq \rho \triangleq {\rm Diam} D < \infty$.
If we simply take $ \lambda_{k} =\lambda \delta_{l}$ with $ \lambda \in (0, \bar \lambda]$ and
$$
 \bar \lambda =2(1-\beta)/(L\rho^{2}),
$$
 then
\begin{equation} \label{eq:4.8}
 f (x^{k}+\lambda_{k} d^{k})
 \leq f (x^{k})+\beta \lambda_{k} \langle f'(x^{k}),d^{k} \rangle,
\end{equation}
as desired; cf. (\ref{eq:4.2}). This means that we can drop the line-search procedure in
Step 2. We call this modification (CGMIL). Obviously, the assertions of
Proposition \ref{pro:4.1} and Theorem
\ref{thm:4.1} remain true for this version.

As (CGMIL) has a two-level structure with each stage containing a
finite number of inner iterations, it is more suitable to derive its
complexity estimate, which gives the total amount of work of the
method. Given a starting point $ w^{0}$ and a
number $\varepsilon > 0$, we define the complexity of the method, denoted
by $N(\varepsilon )$,  as the total number of inner iterations at
$p(\varepsilon )$ stages such that $p(\varepsilon )$ is the maximal number $p$
with $\Delta(w^{p}) \geq \varepsilon$, hence,
\begin{equation} \label{eq:4.9}
 N (\varepsilon ) \leq  \sum ^{p(\varepsilon  )} _{p=1} N_{p},
\end{equation}
where $N_{p}$ denotes the total number of iterations at  stage
$p$. We have to estimate the right-hand side of (\ref{eq:4.9}).

%4===================thm 4.2=========================================

\begin{theorem} \label{thm:4.2}
Let a sequence $\{w^{l}\}$ be generated by  (CGMIL) with the rule:
\begin{equation} \label{eq:4.10}
 \delta _{p} = \nu ^{p}\delta_{0}, \  p=0,1,\ldots;
  \quad \nu \in (0,1), \delta_{0}>0.
\end{equation}
Suppose that the assumptions in {\bf (H)} be fulfilled and also that
the function $f$ is convex and its gradient is Lipschitz continuous  with constant
$L$.  Then the method has the complexity
estimate
$$
N (\varepsilon ) \leq C_{1} \nu ((\delta_{0}/\varepsilon)-1)/(1-\nu),
$$
where $C_{1}=\rho^{2}L/(2\beta(1-\beta)\delta_{0})$.
\end{theorem}
{\bf Proof.}
From (\ref{eq:4.8}) and (\ref{eq:4.1}) we have
$$
 f (x^{k+1}) \leq f (x^{k})-\beta \lambda_{k} \delta _{p}=  f (x^{k})-\beta \bar \lambda \delta^{2} _{p},
$$
at any fixed stage $p$. It follows from the definition of $\bar \lambda$ that
\begin{equation} \label{eq:4.11}
 N_{p} \leq (f (w^{p-1})-f^{*})/(\beta \bar \lambda \delta^{2} _{p})
     \leq \rho^{2}L \Delta(w^{p-1})/(2\beta(1-\beta)\delta^{2} _{p}).
\end{equation}
By the convexity of $f$, for some $ x^{*} \in D^{*}$ we have
$$
\Delta(w^{p}) = f  (w^{p}) - f(x^{*}) \leq \langle f' (w^{p}),w^{p}-x^{*} \rangle  \leq \delta_{p}.
$$
Using this estimate in (\ref{eq:4.11}) gives
$$
N_{p} \leq \rho^{2}L \delta_{p-1}/(2\beta(1-\beta)\delta^{2} _{p}).
$$
From (\ref{eq:4.10}) it now follows that
$$
N_{p} \leq \rho^{2}L \nu ^{-p}/(2\beta(1-\beta)\delta_{0}\nu)=C_{1}\nu ^{-p-1}.
$$
On the other side, since  $\varepsilon \leq \Delta(w^{p}) \leq  \delta_{p}= \delta_{0}\nu ^{p}$, we have
$$
\nu^{-p(\varepsilon  )} \leq \delta_{0}/\varepsilon.
$$
Combining both the inequalities in  (\ref{eq:4.9}), we obtain
\begin{eqnarray*}
N (\varepsilon ) && \leq C_{1}\sum ^{p(\varepsilon  )}
_{p=1} \nu^{-p-1}  = C_{1} \nu (\nu^{-p(\varepsilon  )}-1)/(1-\nu)  \\
&& \leq C_{1} \nu (( \delta_{0}/\varepsilon)-1)/(1-\nu).
\end{eqnarray*}
\QED

We observe that the above estimate is the same as in (\ref{eq:4.7}),
which corresponds to the the usual (CGM) under the same assumptions.

%55555555555555555555555555555555555555555555555555555555555555555555555555555

\section{A parametric inexact method without line-search} \label{sc:5}

In this section, we describe the combined method, which involves both
the inexact solution of the auxiliary direction finding subproblem (\ref{eq:1.2})
due to special parametric threshold control of the descent property and
the simple adaptive step-size rule without line-search.

%===========================Method (CGMIS)==================================
\medskip
\noindent {\bf Method (CGMIS).} \\
 {\em Initialization:} Choose a point $w^{0} \in D$, numbers $\beta \in (0,1)$,
$\theta  \in (0,1)$, and a positive sequence $\{\delta _{p}\} \to 0$. Set $p=1$.

{\em Step 0:} Choose a sequence $\{\tau _{l,p}\} \to 0$, $\tau _{0,p} \in (0, 1)$.
Set  $k=0$, $l=0$, $x^{0}=w^{p-1}$, choose a number $\lambda_{0} \in (0, \tau _{0,p}]$.

{\em Step 1:}  Find a point $z^{k}\in D$ such that
\begin{equation} \label{eq:5.1}
\langle f'(x^{k}), x^{k}-z^{k} \rangle  \geq \delta _{p}.
\end{equation}
If  $\mu(x^{k})<\delta _{p} $, set $w^{p}=x^{k}$, $p=p+1$ and go to Step 0. {\em (Restart)}

{\em Step 2:}  Set $d^{k}= z^{k}-x^{k}$, $x^{k+1}=x^{k}+\lambda_{k}d^{k}$. If
\begin{equation} \label{eq:5.2}
  f(x^{k+1}) \leq f(x^{k})+\beta \lambda_{k}\langle f'(x^{k}), d^{k}\rangle,
\end{equation}
take $\lambda_{k+1} \in [\lambda_{k}, \tau _{l,p}]$. Otherwise
set $\lambda'_{k+1} = \min \{\lambda_{k}, \tau _{l+1,p}\}$,
$l=l+1$ and take $\lambda_{k+1} \in (0,\lambda'_{k+1}]$.
Set $k=k+1$ and go to Step 1.
\medskip

Again, each outer iteration (stage) $p$ contains some number of inner iterations in $k$
with the fixed tolerance $\delta _{p}$. Completing each stage,
that is marked as restart, leads to decrease of its value.
Note that the choice of the parameters $\{\tau _{l,p}\}$ can be in principle independent
for each stage $p$.

By (\ref{eq:5.1}), we again have
$$
\langle f'(x^{k}),d^{k}\rangle \leq -\delta_{p}<0
$$
in (\ref{eq:5.2}). It follows that
 \begin{equation} \label{eq:5.3}
 f (x^{k+1})-f (x^{k}) \leq \beta \lambda_{k} \langle f'(x^{k}),d^{k}\rangle \leq -\beta
 \lambda_{k}\delta_{p}.
\end{equation}

We show that each stage is well defined.

%=========================pro:5.1==================================================

\begin{proposition} \label{pro:5.1} Let the assumptions in {\bf (H)} be fulfilled.
Then the number of iterations at each stage $p$ is finite.
\end{proposition}
{\bf Proof.}
Fix any $p$ and suppose that the sequence $\{x^{k}\}$ is infinite.
Then the number of changes of index $l$ is also infinite.
In fact, otherwise we have $\lambda _{k} \geq \bar \lambda>0$ for $k$ large enough, hence
(\ref{eq:5.3}) gives
$$
f^{*}\leq  f(x^{k+t}) \leq f(x^{k})-t \beta \bar \lambda \delta_{p} \ \to -\infty \ \mbox{as} \  t \to \infty,
$$
for $k$ large enough, which is a contradiction.
Therefore, there exists an infinite subsequence of indices $\{k_{l}\}$
 such that
$$
f (x^{k_{l}}+\lambda _{k_{l}}d^{k_{l}})-f (x^{k_{l}})=f
(x^{k_{l}+1})-f (x^{k_{l}}) > \beta \lambda _{k_{l}}\langle f
'(x^{k_{l}}), d^{k_{l}}\rangle,
$$
or equivalently,
\begin{equation} \label{eq:5.4}
\frac{f (x^{k_{l}}+\lambda _{k_{l}} d^{k_{l}})-f
(x^{k_{l}})}{\lambda _{k_{l}}}
> \beta \langle f '(x^{k_{l}}), d^{k_{l}}\rangle,
\end{equation}
where $d^{k_{l}}=z^{k_{l}}-x^{k_{l}}$. Besides, it holds that
$$
 \lambda_{k_{l}} \in (0,\tau_{l,p}], \ \lambda_{k_{l}+1} \in (0,\tau_{l+1,p}],
$$
where
$$
\lim \limits_{s\rightarrow \infty}\tau_{l,p}=0.
$$
Both the sequences $\{x^{k}\}$ and $\{z^{k}\}$ belong to the bounded
set $D$ and hence have limit points. Without loss of generality, we can
suppose that the subsequence $\{x^{k_{l}}\}$ converges to a point
$\bar x$ and the corresponding subsequence $\{z^{k_{l}}\}$ converges to a point
$\bar z$. Due to (\ref{eq:5.1}) we have
\begin{equation} \label{eq:5.5}
\langle f'(\bar x),\bar z -\bar x \rangle= \lim \limits_{l\rightarrow \infty}
\langle f'(x^{k_{l}}), z^{k_{l}}-x^{k_{l}} \rangle \leq -\delta_{p}.
\end{equation}
At the same time, taking the limit $l\rightarrow \infty$ in (\ref{eq:5.4}), we obtain
\begin{eqnarray*}
\langle f'(\bar x),\bar z -\bar x \rangle &=& \lim \limits_{s\rightarrow \infty
} \left\{\lambda_{k_{l}}^{-1 }(f(x^{k_{l}}+(\lambda_{k_{l}}/\theta) d^{k_{l}}) - f
(x^{k_{l}}))\right\} \\
&\geq&  \beta \langle f'(\bar x),\bar z -\bar x
\rangle,
\end{eqnarray*}
i.e.,  $   (1-\beta ) \langle f'(\bar x),\bar z -\bar x \rangle \geq 0$,
which is a contradiction with (\ref{eq:5.5}). \QED

We are ready to prove convergence of the whole method.
Although it is similar to Theorem \ref{thm:4.1}, we give the full proof for more clarity.

%4===================thm 5.1=========================================

\begin{theorem} \label{thm:5.1}
Let the assumptions in {\bf (H)} be fulfilled. Then:

(i)  The number of changes of index $k$ at each stage $p$  is finite.

(ii) The sequence $\{w^{p}\}$ generated by method (CGMIS) has limit points, all
these limit points belong to $D^{0}$.

(iii)  If $f$ is pseudo-convex, then all the limit points of the sequence $\{w^{p}\}$
belong to the set $D^{*}$, besides, (\ref{eq:4.4a}) holds.
\end{theorem}
{\bf Proof.} Assertion (i) has been obtained in Proposition \ref{pro:5.1}.
By construction, the sequence $\{w^{p}\}$ is bounded, hence it has limit points.
For each $p$ and any point $u^{p}\in Z(w^{p})$ it holds that
\begin{equation} \label{eq:5.6}
\langle f'(w^{p}),w^{p}-u^{p}\rangle \leq \delta_{p}.
\end{equation}
Fix this sequence $\{u^{p}\}$. It is also bounded and must have limit points.
Take an arbitrary limit point $\bar w$ of $\{w^{p}\}$. Then,
without loss of generality we can suppose that
$$
\bar u =\lim \limits_{t\rightarrow \infty }u^{p_{t}} \ \mbox{and} \ \bar w =\lim \limits_{t\rightarrow \infty }w^{p_{t}},
$$
for some subsequences $\{u^{p_{t}}\}$ and $\{w^{p_{t}}\}$.
Taking the limit $t\rightarrow \infty$ in (\ref{eq:5.6}) with $p=p_{t}$, we obtain
$$
\langle f'(\bar w),  \bar w-\bar u \rangle =\lim \limits_{t\rightarrow \infty }
\langle f'(w^{p_{t}}),w^{p_{t}}-u^{p_{t}}\rangle \leq 0.
$$
By Lemma \ref{lm:2.2} we have $\bar u \in Z(\bar w)$, hence
$\langle f'(\bar w),  \bar w-\bar u \rangle =0$ and $\bar w \in D^{0}$
due to Lemma \ref{lm:2.3}.  This means that
all the limit points of $\{u^{p}\}$ belong to $D^{0}$.
This gives assertion (ii).  If $f$ is pseudo-convex, then $D^{0}=D^{*}$
due to Lemma \ref{lm:2.1}, which gives (\ref{eq:4.4a}). We conclude that
assertion (iii) is true.
\QED

Comparing Theorems \ref{thm:3.1} and \ref{thm:5.1}, we observe that
the joint modifications enable us to attain strengthened convergence
properties for (CGMIS) with respect to (CGMS) in the non-convex case.

In this paper, we describe modifications for the basic
conditional gradient method. Obviously, the same modifications
can be applied to most of the gradient type smooth optimization methods.

%66666666666666666666666666666666666666666666666666666666666666666666666666666666666666666666666666

\section{Computational experiments}\label{sc:6}

In order to check the performance of the proposed methods we carried
out computational experiments.  We compared  (CGM), (CGMS), (CGMI),
and (CGMIS) with respect to (\ref{eq:1.1}) for different dimensionality.
They were implemented in Delphi with double precision
arithmetic. The main goal was to compare the number of iterations (it),
the total number of calculations of the goal function value (kf),
and the total number of calculations of  partial derivatives of $f$
(kg) for attaining the same
accuracy $\varepsilon=0.1$ with respect to the gap function
$\mu(x)$. We chose the rule $\delta_{p+1}=\nu \delta_{p}$
 with $\nu = 0.5$ for (CGMI) and (CGMIS). For
 (CGMS) and (CGMIS), we simply set $\lambda_{k+1} =\lambda_{k}$
 if  (\ref{eq:3.1}) (respectively, (\ref{eq:5.2})) holds, and $\lambda_{k+1} =\sigma \lambda_{k}$
 with $\sigma = 0.9$  otherwise.  Next, in the case of restart in
 (CGMIS) we took  $\lambda_{0} =\tau _{0,p}=\lambda_{k}/\sigma$, where
 $\lambda_{k}$ was the current step-size from the previous stage.
 We set $\beta =\theta =0.5$ for all the methods.

We took the simplex as the feasible set for all the test problems, i.e.,
$$
D=\left\{ x\in \mathbb{R}^{m}_{+} \ \vrule \ \sum \limits_{i=1}^{m} x_{i}=b \right\}.
$$
We set $b=10$ and took the same starting point  $x'=(b/m)e$ for all the methods.
For (CGMI) and (CGMIS), we apply the cyclic selection of indices.
In all the series, we took the convex cost functions.

In the first series, we chose $f(x)=\varphi_{1}(x)$ where
\begin{equation} \label{eq:6.2}
\varphi_{1} (x)= 0.5 \langle Px,x \rangle,
\end{equation}
the elements of the matrix $P$ are defined by
\begin{equation} \label{eq:6.3}
p_{ij}= \left\{ {
\begin{array}{rl}
\displaystyle
\sin (i) \cos (j) \quad & \mbox{if} \ i<j, \\
\sin (j) \cos (i) \quad & \mbox{if} \ i>j, \\
\sum \limits_{s \neq i} | p_{is}| +1 \quad & \mbox{if} \ i=j.
\end{array}
} \right.
\end{equation}
The results  are given in Table \ref{tbl:1}.
\begin{table}
\caption{Quadratic cost function $\varphi_{1}$} \label{tbl:1}
\begin{center}
\begin{tabular}{|r|rrr|rrr|}
\hline
  & {} & (CGM) & {} & {} & (CGMS) & {}  \\
\hline
  $n$                &  it & kf & kg  & it & kf & kg   \\
\hline
    $5$            & 202 & 2098 &1010  &  65 & 65 & 325  \\
\hline
     $10$           & 713 & 8256 &7130  &  87 & 87 & 870  \\
\hline
    $20$           & 503 & 6500 & 10060  &  833 & 833 & 16660 \\
\hline
$50$           & 1729 & 24624 & 86450  &  2155 & 2155 & 107750  \\
\hline
$100$         & 2540 & 38454 & 254000  &  5430 & 5430 & 543000  \\
\hline
& {} & (CGMI) & {}  & {} & (CGMIS) & {}   \\
\hline
   $n$              &  it & kf & kg  & it & kf & kg  \\
\hline
    $5$            & 199 & 2072& 472 & 71 & 71 & 192 \\
\hline
     $10$           & 743 & 8690 &4098 & 743 & 743 & 4060 \\
\hline
    $20$            & 583 & 7687 & 5953 & 124 & 124 & 1296 \\
\hline
$50$           & 2037 & 29558 & 50532 & 2214 & 2214 & 53868 \\
\hline
$100$          & 2888 & 44296 & 14864 & 6417 & 6417 & 21764 \\
\hline
\end{tabular}
\end{center}
\end{table}
In the second series, we took the cost function
$$
f(x)=\varphi_{1}(x)+\varphi_{2}(x)
$$
where the function $\varphi_{1}$ was defined as
in (\ref{eq:6.2})--(\ref{eq:6.3}) and the
function $\varphi_{2}$ was defined by the formula
\begin{equation} \label{eq:6.4}
\varphi_{2}(x)=1/(\langle c,x\rangle+d),
\end{equation}
where the elements of the vector $c$ are defined by
$$
c_{i}= 2+\sin(i) \ \mbox{ for } \ i=1,\ldots, m,
$$
and $d=5$.
The results  are given in Table \ref{tbl:2}.
\begin{table}
\caption{Convex cost function $\varphi_{1}+\varphi_{2}$} \label{tbl:2}
\begin{center}
\begin{tabular}{|r|rrr|rrr|}
\hline
 & {} & (CGM) & {} & {} & (CGMS) & {}   \\
\hline
    $n$            &  it & kf & kg  & it & kf & kg   \\
\hline
    $5$            & 203 & 2116 & 1015  &  65 & 65 & 325  \\
\hline
     $10$           & 705 & 8155 & 7050  &  87 & 87 & 870  \\
\hline
    $20$           & 491 & 6329 & 9820  &  833 & 833 & 16660 \\
\hline
$50$           & 1760 & 25105 & 88000  &  2155 & 2155 & 107750  \\
\hline
$100$         & 2594 & 39338 & 259400  &  5496 & 5496 & 549600  \\
\hline
& {} & (CGMI) & {}  & {} & (CGMIS) & {}   \\
\hline
    $n$              &  it & kf & kg  & it & kf & kg  \\
\hline
    $5$            & 209 & 2192 & 517 & 71 & 71 & 192 \\
\hline
     $10$           & 731 & 8528 & 4008 & 706 & 706 & 3866 \\
\hline
    $20$            & 547 & 7150 & 5660 & 123 & 123 & 1291 \\
\hline
$50$           & 2026 & 29359 & 50551 & 2070 & 2070 & 50506 \\
\hline
$100$          & 2921 & 44797 & 16343 & 6321 & 6321 & 17517 \\
\hline
\end{tabular}
\end{center}
\end{table}

In the third series, we chose $f(x)=\varphi_{3}(x)$ where
\begin{equation} \label{eq:6.5}
\varphi_{3} (x)= 0.5 \| Px-q \|^{2},
\end{equation}
the elements of the $m \times n $ matrix $P$ are defined by
\begin{equation} \label{eq:6.6}
p_{ij}= \left\{ {
\begin{array}{ll}
\displaystyle
\tilde p_{ij} \quad & \mbox{if} \ i \neq j, \\
\tilde p_{ij}+2 \quad & \mbox{if}  \ i=j;
\end{array}
} \right.
\end{equation}
where
\begin{equation} \label{eq:6.7}
 \tilde p_{ij}= \ln(1+i/j)\sin (i/j)/ (i+j), \ i=1, \dots, m, j=1, \dots, n;
\end{equation}
and
\begin{equation} \label{eq:6.8}
q_{i}= b \sum  \limits_{j=1}^{n}  p_{ij}, \ i=1, \dots, m.
\end{equation}
 The results  are given in Table \ref{tbl:3}.
\begin{table}
\caption{Quadratic cost function $\varphi_{3}$} \label{tbl:3}
\begin{center}
\begin{tabular}{|rr|rrr|rrr|}
\hline
& {}  & {} & (CGM) & {} & {} & (CGMS) & {}  \\
\hline
 $m$ & $n$                &  it & kf & kg  & it & kf & kg   \\
\hline
  $2$ &  $5$            & 2362 & 25584 & 11810  &  68 & 68 & 340  \\
\hline
  $5$ &   $10$           & 2998 & 33752 & 29980  &  88 & 88 & 880  \\
\hline
  $10$ &  $20$           & 4076 & 43381 & 81520  &  77 & 77 & 1540 \\
\hline
$25$ & $50$           & 219 & 2438 & 10950  &  497 & 497 & 24850  \\
\hline
$50$ & $100$         & 4197 & 48491 & 419700  &  93 & 93 & 9300  \\
\hline
& {} & {} & (CGMI) & {}  & {} & (CGMIS) & {}   \\
\hline
 $m$ &  $n$              &  it & kf & kg  & it & kf & kg  \\
\hline
  $2$ &  $5$            & 2412 & 26273 & 3617 & 54 & 54 & 103 \\
\hline
   $5$ &  $10$           & 3328 & 38601 & 9379 & 107 & 107 & 255 \\
\hline
 $10$ &   $20$            & 4296 & 50567 & 18508 & 93 & 93 & 393 \\
\hline
$25$ & $50$           & 3924  & 45769 & 35894 & 837 & 837 & 7343 \\
\hline
$50$ & $100$          & 4176 & 48969 & 5385 & 116 & 116 & 1774 \\
\hline
\end{tabular}
\end{center}
\end{table}
In the fourth series, we took the cost function
$$
f(x)=\varphi_{3}(x)+\varphi_{2}(x)
$$
where the function $\varphi_{3}$ was defined as
in (\ref{eq:6.5})--(\ref{eq:6.8}) and the function $\varphi_{2}$ was defined as
in (\ref{eq:6.4}).
The results  are given in Table \ref{tbl:4}.
\begin{table}
\caption{Convex cost function $\varphi_{3}+\varphi_{2}$} \label{tbl:4}
\begin{center}
\begin{tabular}{|rr|rrr|rrr|}
\hline
& {}  & {} & (CGM) & {} & {} & (CGMS) & {}  \\
\hline
 $m$ & $n$                &  it & kf & kg  & it & kf & kg   \\
\hline
  $2$ &  $5$            & 2365 & 25633 & 11825  &  68 & 68 & 340  \\
\hline
  $5$ &   $10$           & 2926 & 32957 & 29260  &  93 & 93 & 930  \\
\hline
  $10$ &  $20$           & 3972 & 46047 & 79440  &  83 & 83 & 1660 \\
\hline
$25$ & $50$           & 229 & 2564 & 11450  &  497 & 497 & 24850  \\
\hline
$50$ & $100$         & 4192 & 48894 & 419200  &  92 & 92 & 9200  \\
\hline
& {} & {} & (CGMI) & {}  & {} & (CGMIS) & {}   \\
\hline
 $m$ &  $n$              &  it & kf & kg  & it & kf & kg  \\
\hline
  $2$ &  $5$            & 2386 & 25956 & 3577 & 54 & 54 & 103 \\
\hline
   $5$ &  $10$           & 3276 & 37911 & 9271 & 107 & 107 & 255 \\
\hline
 $10$ &   $20$            & 4543 & 53856 & 19532 & 92 & 92 & 392 \\
\hline
$25$ & $50$           & 3806  & 44284 & 34835 & 867 & 867 & 7615 \\
\hline
$50$ & $100$          & 4304 & 50627 & 7408 & 115 & 115 & 1763 \\
\hline
\end{tabular}
\end{center}
\end{table}

In almost all the cases, (CGMS) and (CGMIS), which do not use line-search, showed
rather rapid convergence, they outperformed (CGM) and (CGMI), respectively,
 in the total number of  goal function calculations. Similarly, the inexact versions
 (CGMI) and (CGMIS) showed essential reduction of the total number of  partial derivatives
 calculations in comparison with (CGM) and (CGMS), respectively.

%%%%%%%%%%%%%%%%%%%%%%%%%%%%%%%%%%%%%%%%%%%%%%%%%%%%%%%%%%%%%%%%%%%%%%%%%%%%%%%%%%%%%

\section*{Acknowledgement}

The results of this work were obtained within the state assignment of the
Ministry of Science and Education of Russia, project No. 1.460.2016/1.4.
In this work, the author was also supported by Russian Foundation for Basic Research, project No.
13-01-00368-a.

%##########################################################################


\begin{thebibliography}{99}



\bibitem{FW56}
M.~Frank and P.~Wolfe, {\em An algorithm for quadratic programming}, Nav.
Res. Logist. Quart., vol.~3 (1956), pp.95--110.

\bibitem{LP66}
E.S.~Levitin and B.T.~Polyak, {\em Constrained minimization
methods}, USSR Comp. Maths. Math. Phys., vol.~6 (1966), pp.1--50.

\bibitem{DR68}
 V.F.~Dem'yanov and A.M.~Rubinov, {\em Approximate Methods for Solving
Extremum Problems}, Leningrad Univ. Press, Leningrad, 1968. (Engl.
transl. in Elsevier, Amsterdam, 1970)

\bibitem{PD78}
B.N.~Pshenichnyi and Yu.M.~Danilin,
{\em Numerical Methods in Extremal Problems}, MIR, Moscow, 1978.

\bibitem{Dun80}
J.C.~Dunn,  {\em Convergence rates for conditional gradient sequences generated by implicit step length
rules}, SIAM J. Control Optim., vol.~18 (1980), pp.473--487.

\bibitem{Cla10}
K.L.~Clarkson,  {\em Coresets, sparse greedy approximation, and the
Frank-Wolfe algorithm}, ACM Trans. on Algor., vol.~6 (2010), Art. No. 63,
pp.1--30.

\bibitem{Jag13}
M.~Jaggi, {\em  Revisiting Frank-Wolfe: Projection-free sparse convex
optimization},  Proc. of the 30th International Conference on Machine
Learning (ICML-13),  2013, pp.427--435.

\bibitem{GM86} J.~Guelat and P.~Marcotte, {\em Some comments on Wolfe's \lq\lq away step"}, Math.
Progr., vol.~35 (1986), pp.110--119.

\bibitem{BT04}
A. Beck and M. Teboulle, {\em A conditional gradient method with linear rate
of convergence for solving convex linear systems}, Math. Meth. Oper.
Res.,  vol.~59 (2004), pp.235--247.

\bibitem{FG16}
R.M.~Freund and P.~Grigas, {\em  New analysis and results for the
Frank-Wolfe method}, Mathem. Progr., vol.~155 (2016),  pp.199--230.

\bibitem{Kon17}
I.V.~Konnov, {\em The method of pairwise variations with tolerances for
linearly constrained optimization problems},
J. Nonlin. Variat. Anal., vol.~1 (2017), pp.25--41.

\bibitem{Kon18}
I.V.~Konnov, {\em Conditional gradient method without line-search},
Russ. Mathem. (Izv. VUZ), vol.~62 (2018), pp.82--85.

\bibitem{Kon13}
I.V.~Konnov, {\em Nonlinear Optimization  and Variational
Inequalities}, Kazan Univ. Press, Kazan, 2013. [In Russian]

\bibitem{DH78}
J.C. Dunn and S. Harshbarger, {\em Conditional gradient algorithms with open
loop step size rules}, J. Math. Anal. Appl., vol.~62 (1978), pp. 432--444.

\bibitem{Kon16b}
I.V.~Konnov, {\em  Selective bi-coordinate variations for resource
allocation type problems}, Comp. Optim. Appl., vol.~64 (2016), pp.821--842.

\bibitem{Kon17a}
I.V.~Konnov, {\em  An adaptive partial linearization method for optimization
problems on product sets}, J. of Optim. Theory and Appl.,
vol.~175 (2017), pp.478--501.

\end{thebibliography}
\end{document}